\input amstex
\documentstyle{amsppt}
\magnification=\magstep1

\baselineskip 15pt \hsize=12.3 cm \vsize=18.5 cm \pagewidth{12.3
cm} \pageheight{18.5 cm} \loadbold \TagsOnRight \NoBlackBoxes

\def\build#1_#2^#3{\mathrel{\mathop{\kern 0pt#1}\limits_{#2}^{#3}}}

\def\buld#1^#2{\mathrel{\mathop{\kern 0pt#1}\limits^{#2}}}

\font\bforteen=cmbx9 at 13pt \font\bfifteen=cmbx9 at 15pt

\def\picture #1 by #2 (#3){\vbox to #2{\hrule width #1
height 0pt\vfill \special{picture #3}}}

\def\scaledpicture #1 by #2 (#3 scaled #4){{
\dimen0=#1 \dimen1=#2 \divide\dimen0 by 1000 \multiply\dimen0 by
#4 \divide\dimen1 by 1000 \multiply\dimen1 by #4 \picture \dimen0
by \dimen1 (#3 scaled #4)}}

\topmatter
\title\nofrills{\bfifteen
Hitting and return times in ergodic dynamical systems}
\endtitle
\author
N\. Haydn, Y\. Lacroix, \& S. Vaienti
\endauthor

\affil Los Angeles, Toulon, \& Marseille
\endaffil
\leftheadtext{Hitting and return times} \rightheadtext{N. Haydn,
Y. Lacroix \& S. Vaienti}

 \keywords asymptotic
distribution, hitting, return times, Kac
\endkeywords

\abstract\nofrills\centerline{\bf Abstract} Given an ergodic
dynamical system $(X,T,\mu )$, and $U\subset X$ measurable with
$\mu (U)>0$, let $\mu (U)\tau_U(x)$ denote the normalized hitting
time of $x\in X$ to $U$. We prove that given a sequence $(U_n)$
with $\mu (U_n)\to 0$, the distribution function of the normalized
hitting times to $U_n$ converges weakly to some sub-probability
distribution $F$ if and only if the distribution function of the
normalized return time converges weakly to some distribution
function $\tilde F$, and that in the converging case,
$$
F(t)=\int_0^t(1-\tilde F(s))ds,\; t\ge 0.\tag$\diamondsuit$
$$
This in particular characterizes asymptotics for hitting times,
and shows that the asymptotics for return times is exponential if
and only if the one for hitting times is too.
\endabstract

\nologo \TagsOnLeft
\endtopmatter

\document

\medskip

\noindent{\bforteen 1. Introduction.}

\medskip

In the recent years there has been an interest in the statistics
of hitting and return times. Typically a neighbourhood of a point
is considered which can be either a metric ball or a `cylinder
set' associated with a measurable partition. In accordance with a
theorem due to Kac one then looks at the return times which are
normalized by the measure of the return set. A number of recent
papers (e.g. [A1], [A2], [AG], [BV1], [BV2], [Co], [CG1], [CO],
[GS], [H], [H1], [H2], [HSV], [HV], [P], [PI]) have provided
conditions under which this distribution converges to the
exponential distribution if the set is shrunk so that its measure
converges to zero. On a different note, Lacroix and Kupsa have
shown that with a suitable choice of return set one can realize
any arbitrarily chosen limiting return time distribution [L] and
hitting time distribution [K-L] within some class (see Theorem 1).

The purpose of this note is to show that limiting distributions
for hitting and return times are intimately linked by the
transformation $(\diamondsuit)$.

Let $(X,\Cal B,\mu)$ be a probability space, $T\colon X\to X$ a
measurable transformation that preserves $\mu$, i.e.\ $T^*\mu
=\mu$. We also assume the dynamical system $(X,T,\mu )$ to be
ergodic.

For (measurable) $U\subset X$ with $\mu (U)>0$ we define the {\bf
return/hitting time function} $\tau_U$ by
$$
\tau_U(x)=\inf\{k\ge 1\colon T^kx\in U\}.
$$
For $x\in U$, $\tau_U(x)$ denotes the {\bf return time}. On the
other hand if we refer to $\tau_U(x)$ as a function on all of $X$
then we call it the {\bf hitting time function}. Poincar\'e's
recurrence theorem [K, Theorem 1'], together with ergodicity, then
shows that $\tau_U$ is a.s. finite. We also have Kac's theorem [K,
Theorem 2'] according to which
$$
\int_U \tau_U(x)\,d\mu(x)=\sum_{k=1}^\infty
k\mu(U\cap\{\tau_U=k\})=1.
$$

Finer statistical properties of the variable $\mu(U)\tau_U$ have
been investigated, in a rather large number of recent papers,
where particular attention was given to the study of weak
convergence of $\mu(U_n)\tau_{U_n}$ as $\mu(U_n)\to 0$. See [AG]
for a recent survey in the mixing case.

We say a sequence of distribution functions $F_n$, $n=1,2,\dots$,
{\bf converges weakly} to a function $F$ (which might not be a
distribution itself) if $F$ is (non-strictly) increasing right
continuous and satisfies $\lim_{n\rightarrow\infty}F_n(t)=F(t)$ at
every point $t$ of continuity of $F$. We will write
$F_n\Rightarrow F$ if $F_n$ converges weakly to $F$.

\medskip

Given a $U\subset X$ measurable with $\mu (U)>0$, we define
$$
\matrix \tilde F_{U}(t)=\frac1{\mu
(U)}\mu\left(U\cap\left\{\tau_{U}\mu (U)\le
t\right\}\right)\\
F_{U}(t)=\mu(\{\mu(U)\tau_U\le t\}).
\endmatrix
$$

Define
$$
\left\{
\matrix
\Cal F=
\left\{F:\Bbb R\to [0,1],\; F\equiv 0\text{
on $]-\infty ,0]$, $F$ (non-strictly) increasing, }\right.\hfill\\
\left.\text{\qquad\qquad\qquad\qquad\qquad continuous, concave on
$[0,+\infty [$,
$F(t)\le t$ for $t\ge 0$}\right\};\hfill\\
\tilde\Cal F=
\left\{\tilde F:\Bbb R\to [0,1],\text{ $\tilde F$
(non strictly) increasing, right-continuous,} \right.\hfill\\
\left.\qquad\qquad\qquad\qquad\tilde F\equiv 0\text{ on }]-\infty
,0],
\;\int_0^{+\infty}(1-\tilde F(s))ds\le 1\right\}.\hfill\\
\endmatrix
\right.
$$

These functional classes appear in the following ($U_n$ is always
assumed to be measurable):

\proclaim{Theorem 1} Let $(X,T,\mu )$ be an ergodic and aperiodic
dynamical system. Then:

(a) {\rm [L]} for any $\tilde F\in\tilde \Cal F$ there exists a
sequence $\{U_n\subset X: n=1,2,\dots\}$ such that $\mu(U_n)\to 0$
and $\tilde F_{U_n}\Rightarrow\tilde F$.

(b) {\rm [L]} if $\mu (U_n)\to 0$ and $\tilde F_{U_n}\Rightarrow
\tilde{F}$ then $\tilde F\in\tilde\Cal F$.

(c) {\rm [K-L]} for any $F\in\Cal F$, there exists $\{U_n\subset
X: n=1,2,\dots\}$ such that $\mu(U_n)\to 0$ and
$F_{U_n}\Rightarrow F$.

(d) {\rm [K-L]} if $\mu (U_n)\to 0$ and $F_{U_n}\Rightarrow {F}$
then $F\in\Cal F$.
\endproclaim

\medskip

In this note we prove the following rather unexpected and
surprisingly unknown result:

\proclaim{Main Theorem} Let $(X,T,\mu )$ be ergodic, and
$\{U_n\subset X: n\ge1\}$ a sequence of positive measure
measurable subsets. Then the sequence of functions $\tilde
F_{U_n}$ converges weakly if and only if the functions $F_{U_n}$
converge weakly.

Moreover, if the convergence holds, then
$$
F(t)=\int_0^t(1-\tilde F(s))ds,\; t\ge 0,\tag $\diamondsuit$
$$
where $\tilde F$ and $F$ are the corresponding limiting
(sub-probability) distributions.
\endproclaim

The only previous result in this direction was obtained in [HSV]
where it is shown that $\tilde F_{U_n}\rightarrow\tilde F$ and
$\tilde F(t)=1-e^{-t}$ for $t\ge 0$ if and only if $F_{U_n}-\tilde
F_{U_n}\to 0$ in the supremum norm on the real line. Our Main
Theorem shows that the exponential distribution is the only
distribution which can be asymptotic to both return and hitting
times, as it is clearly the only fixed point of $(\diamondsuit)$.
The equation $(\diamondsuit)$ also gives an equivalence between
parts (a) and (b) of Theorem 1. We state a Corollary is this
direction :

\proclaim{Corollary 2} (i) The asymptotic distribution for hitting
times, if it exists, is positive exponential with parameter $1$ if
and only if the one for return times is too.

(ii) Parts (a) and (b) of Theorem 1 are equivalent.
\endproclaim

\medskip

Before we prove the Main Theorem in Section 3, and Corollary 2, we
present an application to irrational rotations on the torus.

\medskip

\noindent{\bforteen 2. Homeomorphisms of the circle. }

\medskip

In the paper [CF] (see also [C] for further developments), Coelho
and De Faria  were able to characterize the limiting laws for the
distribution of hitting time for orientation preserving
homeomorphisms of the circle without periodic points, provided the
set $U$ is chosen in a descending chain of renormalisation
intervals. We could therefore apply our theorem to compute the
corresponding distribution for the first return time.

Let $f:S^1\rightarrow S^1 $ be an orientation preserving
homeomorphism of the circle without periodic points and let
$\alpha=\alpha(f)\in[0,1)$ be its irrational rotation number; we
consider here its continued fraction expansion: $\alpha=[a_0,a_1,
\ldots, a_n, \ldots]$.  It is well known that $a_n=[H^n(\alpha)]$,
where $[.]$  denotes now the integer part and $H:
[0,1[\rightarrow[0,1[$ is the Gauss transformation defined by
$H(0)=0$ and $H(\alpha)=\{1/\alpha\}$ for $\alpha \neq 0$,
denoting $\{.\}$ the fractional part. The truncated expansion of
order $n$ of $\alpha$ is given by $p_n/q_n= [a_0,a_1, \ldots,
a_{n-1}]$, where $p_n$ and $q_n$ verify the recursive relations
$$
\matrix
q_{k+1} \!\! & = & \!\! a_k q_k+q_{k-1}\hfill\\
p_{k+1} \!\! & = & \!\! a_k p_k+p_{k-1}\hfill\\
\endmatrix
$$
with $p_0 = 1$, $p_1 = 1$ and $q_0 = 0$, $q_1 = 1$.  Then, for any
number $\beta \in [0, 1]$ and by setting $b_j = [H^j(\beta)]$ for
$j\ge 0$, we can construct the following quantities for $n\ge 1$
(natural extension of $H$)
$$
  \Gamma^n(\alpha, \beta) = (H^n(\alpha),\, [a_{n-1}, a_{n-2}, \ldots,
    a_0, b_0, b_1, \ldots]).
$$
Notice that the convergent subsequences of $ \Gamma^n(\alpha,
\beta)$ for $n\ge 0$ do not depend on $\beta$. The sets $U$, where
the points of the circle enter, are constructed in the following
way. Let us take any point $z$ on the circle and define $J_n$ as
the closed interval of endpoints $f^{q_{n-1}}(z)$ and
$f^{q_n}(z)$: the sequence of sets $U$ shrinking to $z$ is taken
as the descending chain of intervals $J_n$. We call $F_n(t)$ the
distribution function of the hitting time into the sets
$U=\{J_n\}$. Coelho and De Faria proved the following theorem :

\proclaim{Theorem [CF]} For each subsequence $\sigma=\{n_i\}$ of
$\Bbb N$, the corresponding distribution functions $F_{n_{i}}$
converge (pointwise or uniformly) if and only if either

(a) $H^{n_i}(\alpha)\rightarrow 0$, in which case the limit
distribution is the uniform distribution on the unit interval; or

(b) $\Gamma(\alpha,\cdot) \rightarrow (\theta,\omega)$ for some
$\theta>0$ and $\omega<1$, in which case the limit distribution is
the continuous piecewise linear function $F_{\sigma}$ given by
$$
F_{\sigma}(t) = \left\{\matrix
    t  & \text{ if \quad} 0\le t< \frac{(1+\theta)\omega}
    {1+\theta\omega}; \hfill\\
    { \frac{1}{1+\omega}t+\frac{\omega^2(1+\theta)}{(1+\theta\omega)(1+\omega)}}  & \text{  if
    \quad}\frac{(1+\theta)\omega}{1+\theta\omega} \le
    t < \frac{1+\theta}
    {1+\theta\omega}\hfill\\
1 & \text{ if \quad} t \ge \frac{1+\theta}
    {1+\theta\omega}\hfill\\
    \endmatrix\right.
$$
\endproclaim

For the second, more interesting case, we get immediately the
distribution of the first return time by applying our result. It
reads:
$$
  \tilde{F_{\sigma}}(t)  =
  \left\{\matrix
    0  & \text{ if \quad} 0\le t< \frac{(1+\theta)\omega}
    {1+\theta\omega}\hfill\\
\frac{\omega}{1 + \omega}  & \text{  if \quad}
\frac{(1+\theta)\omega}{1+\theta\omega} \le
    t < \frac{1+\theta}
    {1+\theta\omega}\hfill\\
1  & \text{ if \quad} t \ge \frac{1+\theta}
    {1+\theta\omega}
\hfill\\
\endmatrix
\right.
$$

As a concrete example we could take $\alpha=\frac{\sqrt{5}-1}{2}$,
the golden number, which exhibits the continued fraction expansion
$[1,1,1...]$. In this case it is easy to check that
$\Gamma^n(\alpha,.)$ converges, when $n\rightarrow \infty$, to
$(\theta=\alpha, \omega=\alpha)$.

\medskip

\noindent{\bforteen 3. Proof of the Main Theorem.}

\medskip

For $k\ge 1$, let us denote $V_k=\left\{x\in U: \tau_U=k\right\}$.
Then up to a zero measure set $X$ is the disjoint union of the
sets $\bigcup_{j=0}^{k-1}T^{j}V_k$, $k=1,2,\dots$ (this is Kac's
Theorem 1', no invertibility of the transformation is needed).

Let us also introduce $U_k=\{x:\tau_U(x)=k\}$. Then
$$
U_k=\cup_{j=0}^{+\infty}T^jV_{k+j},\text{ and }\mu
(U_k)=\sum_{j\ge k}\mu (V_j).
$$
Also for $t\ge 0$, $F_U(t)=\mu (\{\tau_U\le t/\mu
(U)\})=\sum_{k\le t/\mu (U)}\mu (U_k)$ : notice that $F_U$ is
constant on intervals $[k\mu (U),(k+1)\mu (U)[$, with a jump $\mu
(U_k)$ at $k\mu (U)$.

We next define $\bar F_U$ by
$$
\left\{ \matrix (a)&:&\bar F_U(k\mu (U))=F_U(k\mu(U));\hfill\\
(b)&:&\bar F_U\text{ is linear on }[k\mu (U),(k+1)\mu
(U)].\hfill\\
\endmatrix
\right.
$$
Then because the discontinuity jumps of $F_U$ decrease, the
function $\bar F_U$ is (non-strictly) increasing, concave on the
positive axis, piecewise linear, continuous, and the right hand
side derivative satisfies
$$
{\bar F_U}^{'+}(t)=\frac{\mu (U_{k+1})}{\mu (U)}\text{ for }t\in
[k\mu (U),(k+1)\mu (U)[.
$$
On the other hand, $\tilde F_U(t)=\frac 1{\mu (U)}\sum_{k\le t/\mu
(U)}\mu (V_k)$, whence
$$
\tilde F_U\text{ is constant on }[k\mu (U),(k+1)\mu (U)[ \text{
and has jump }\frac{\mu (V_k)}{\mu (U)}\text{ at }k\mu (U).
$$

Putting the above together yields
$$
{\bar F_U}^{'+} (t)=1-\tilde F_U(t),\;\; t\ge0. \tag $\star$
$$
Notice also that
$$
\parallel F_U-\bar F_U\parallel_{\infty}\le\mu (U)
\tag $\star\star$
$$
since $F_U$ has its discontinuities located at points $\mu (U)$,
$2\mu (U)$, $\ldots$.

\bigskip

We continue with the proof of the Main Theorem:

(I) Let us assume there is a sequence of subsets $U_n\subset X$ so
that $\mu (U_n)\to 0$ and $\tilde F_{U_n}\Rightarrow \tilde F$
(then $\tilde F\in\tilde\Cal F$). Since $\tilde F$ is
(non-strictly) increasing, this implies that $\tilde F_{U_n}\to
\tilde F$ Lebesgue almost surely on $[0,+\infty [$. Whence, for
given $t\ge 0$, by the Lebesgue dominated convergence theorem on
$[0,t]$ ($\tilde F\in [0,1]$), combining with $(\star)$ one has
$$
\bar F_{U_n}(t)=\int_0^t(1-\tilde
F_{U_n}(s))ds\to\int_0^t(1-\tilde F(s))ds =:F(t).
$$
We put $F(t)=0$ for $t< 0$. Since $\tilde F\in\tilde\Cal F$, it
follows that $F\in\Cal F$.

Moreover, by $(\star\star)$, $F_{U_n}(t)\to F(t)$ for all
$t\in\Bbb R$ (the convergence is in fact uniform on compact
subsets of $\Bbb R$ by [R, Theorem 10.8]).

Hence if $\tilde F_{U_n}\Rightarrow\tilde F$, then $(F_{U_n})$
converges weakly to the $F$ associated to $\tilde F$ by formula
$(\diamondsuit )$.

\medskip

Proving the converse for the Main Theorem, we need the following:

\proclaim{Lemma 3} Let $f_n$, $n=1,2,\dots$, be a sequence of
concave functions defined on a non-empty open interval $I$ and
assume that $f_n$ converges pointwise to a limit function $f$.
Then off an at most countable subset of $I$ the sequence of
derivatives $f_n'$ converges pointwise to the derivative $f'$ of
$f$.
\endproclaim

\demo{Proof of Lemma 3} By [R, Theorem 25.3], off an at most
countable subset of $I$, the functions $f_n$, and $f$, are
differentiable, as concave functions.

Next, using the argument for the proof of [R, Theorem 25.7], but
for a fixed $x\in I$ rather than along a sequence of point $x_i$
or points $x_i$ in a closed bounded subset of $I$, the convergence
of the derivatives, when all defined, follows at once.
\qed\enddemo

(II) Let us now assume that $F_{U_n}\Rightarrow F$. Then [KL]
implies that $F\in\Cal F$. Whence by $(\star )$ and
$(\star\star)$, we have, for $t\ge 0$,
$$
\bar F_{U_n}(t)=\int_0^t{\bar F_{U_n}}^{'+}(s)ds=
\int_0^t(1-\tilde F_{U_n}(s))ds\to F(t)\; (=\int_0^t{F'}^+(s)ds).
$$
It now follows from Lemma 3 that off an at most countable subset
$\Omega$ of $]0,+\infty [$, the functions $1-\tilde F_{U_n}(s)$
converge pointwise to $F^{'+}(s)$. We a priori have $F^{'+}(s)$
non-strictly decreasing. We then set $\tilde F(s)=1-F^{'+}(s)$
where $F^{'+}$ is continuous, else we ask $\tilde F$ to be
right-continuous which sets automatically its values at eventual
discontinuity points of $F^{'+}$.

It remains to show that $\tilde F_{U_n}(s)\to \tilde F(s)$ at
points $s$ of continuity of $F^{'+}$. Clearly if $s\notin\Omega$
or $s<0$ there is nothing to do. Else, for any $s_1$ and $s_2$ not
in $\Omega$ such that $s_1<s<s_2$, we have
$$
\tilde F(s_1)\le\liminf_n\tilde F_{U_n}(s)\le\limsup_n\tilde
F_{U_n}(s)\le\tilde F(s_2),
$$
and since $\Omega$ is dense in $[0,+\infty[$, the conclusion
follows. So $\tilde F_{U_n}\Rightarrow \tilde F=1-F^{'+}$, which
ends the proof. \qed

\medskip

\demo{Proof of Corollary 2}

(i) : the proof is standard.

(ii) : we have a map $\tilde F\mapsto F$ defined by
$(\diamondsuit)$ which clearly maps $\tilde\Cal F\to \Cal F$ in a
one to one way.

To prove this map is surjective we may proceed as follows : given
$F\in\Cal F$ we have for any $s\ge 0$, $F(s)=\int_0^sF^{'+}(t)dt$,
where $F^{'+}$ is the righthand side derivative of $F$. Since
$F\in\Cal F$, it follows quite easily that $F^{'+}$ is positive
(non-strictly) decreasing and takes on values in $[0,1]$. By
integrability ($F(s)\to 1$ when $s\to +\infty$), $\lim_{+\infty
}F^{'+} =0$.

We can set $\tilde F$ as to be equal to $1-F^{'+}$ on $]0,+\infty
[$, at continuity points of $F^{'+}$. We else set $\tilde F\equiv
0$ on $]-\infty ,0]$, and we eventually complete the definition of
$\tilde F$ on $]0,+\infty [$ by requiring that it is
right-continuous.

Then one has $\tilde F\mapsto F$ by $(\diamondsuit)$, and $\tilde
F\in\tilde\Cal F$. \qed\enddemo

\bigskip
$$
\matrix
\text{\bf Nicolai Haydn}\hfill&\text{\bf Yves Lacroix}\hfill\\
&&\hfill\\
\text{University of Southern California}\hfill&\text{Universit\'e
du Sud Toulon Var}\hfill\\
\text{Department of
Mathematics}\hfill&\text{ISITV}\hfill\\
\text{Denney Research Building}\hfill&\text{Avenue Georges
Pompidou}\hfill\\
\text{DRB 155}\hfill&\text{BP 56}\hfill\\
\text{Los Angeles}\hfill&\text{83162 La Valette
Cedex}\hfill\\
\text{California 90089-1113}\hfill&\text{France}\hfill&\text{}\hfill\\
\text{USA}\hfill&\text{}\hfill&\text{}\hfill\\
\text{}\hfill&\text{}\hfill&\text{}\hfill\\
\text{nhaydn\@usc.edu}\hfill&\text{ylacroix\@univ-tln.fr}
\hfill\\
&\hfill\\
\text{\bf Sandro Vaienti}\hfill&\hfill\\
&\hfill\\
\text{UMR-6207 Centre de
Physique Th\'eorique, CNRS, }\hfill&\hfill\\
\text{Universit\'es d'Aix-Marseille I, II, et du Sud Toulon Var}\hfill&\hfill\\
\text{FRUMAM (F\'ed\'eration de Recherche des}\hfill&\hfill\\
\text{Unit\'es de Math\'ematiques de Marseille)}\hfill&\hfill\\
\text{CPT}\hfill&\hfill\\
\text{Campus of Luminy}\hfill&\hfill\\
\text{Case 907}\hfill&\hfill\\
\text{13288 Marseille Cedex 9}\hfill&\hfill\\
\text{France}\hfill&\hfill\\
&\hfill\\
\text{vaienti\@cpt.univ-mrs.fr}\hfill&\hfill\\
\endmatrix
$$
\enddocument